\documentclass{article}

\usepackage[linesnumbered,ruled,vlined]{algorithm2e}
\usepackage{amsfonts}
\usepackage{amsmath}
\usepackage{graphicx}
\usepackage{url}
\usepackage{hyperref}
\usepackage[brazil]{babel}

\begin{document}

\author{Samuel L. Azorli \\ FT/Unicamp \\ Limeira-SP \and Luis A. A. Meira \\ FT/Unicamp \\ Limeira-SP}

\title{Modelagem para o Problema de Entrega de Refeições em Rio Claro-SP}

\newcommand{\Sigla}[2]{#1 (#2)}
\newcommand{\textcite}[1]{\cite{#1}}

\newtheorem{definition}{Definição}[section]

\maketitle

\begin{abstract}
Problemas de Roteamento são problemas em que um conjunto de clientes é atendido por um conjunto de veículos. Neste trabalho, modelamos em um mapa 2D um \textit{benchmark} multiobjetivo baseado em um problema de roteamento realístico de entregas de marmitas por motocicletas na cidade de Rio Claro-SP. O mapa gerado de Rio Claro apresenta 1566 ruas com coordenadas extraídas manualmente e modeladas através de cadeias poligonais. Geramos um total de 23 instâncias contendo de 2 a 7 depósitos e até 2000 pontos de entrega. Este trabalho é uma extensão de \cite{2020zeni} no qual os autores modelam o problema de entrega de correspondências por carteiros na cidade de Artur Nogueira.
O trabalho~\cite{2020zeni} possui um depósito, enquanto neste trabalho temos múltiplos depósitos. Em~\cite{2020zeni} foi modelada a cidade de Artur Nogueira, com 537 ruas. Neste trabalho, foi modelado um mapa de Rio Claro, com 1566 ruas. 
As instâncias geradas serão disponibilizadas para que a comunidade científica valide e compare algoritmos de otimização.
\end{abstract}

\begin{abstract}
Vehicle Routing Problems are problems in which a set of customers is visited by a set of vehicles. In this work, we model a multiobjective \textit{benchmark} on a 2D map based on a realistic routing problem of food deliveries by motorcycles in the city of Rio Claro-SP. The generated map of Rio Claro presents 1566 streets with manually extracted coordinates and modeled through polygonal chains. We generate a total of 23 instances containing 2 to 7 deposits and up to 2000 delivery points. This work is an extension of \cite{2020zeni} in which the authors model the problem of mail delivery by postmen in the city of Artur Nogueira. 
The research of~\cite{2020zeni} has one deposit, while in this work we have multiple deposits. In~\cite{2020zeni}, the city of Artur Nogueira was modeled, with 537 streets. In this work, a map of Rio Claro is modeled, with 1566 streets.
The generated instances will be made available for the scientific community to validate and compare optimization algorithms.
\end{abstract}

\section{Introdução}\label{chp:Introducao}

A logística está presente em diversas atividades do cotidiano moderno e é um campo de estudo científico. Problemas de logística podem reduzir custos para empresas, governos e pessoas. A computação trabalha para trazer soluções para problemas do mundo real e auxiliar na tomada de decisão.

Um dos problemas computacionais que envolvem a logística é o famoso problema do caixeiro-viajante (do inglês \Sigla{\textit{Traveling Salesman Problem}}{TSP}). Uma generalização do TSP é o problema de roteamento de veículos (do inglês \Sigla{\textit{Vehicle Routing Problem}}{VRP}). O VRP foi introduzido pelo trabalho \cite{dantzig1959truck}, inicialmente como um problema de abastecimento de combustível por vários caminhões a partir de uma única refinaria. Entretanto o nome VRP surge apenas no trabalho \cite{christofides1976vehicle}, definindo VRP como uma classe genérica de problemas que envolvem a visita de ``clientes''  por ``veículos''.

Desde então surgiram diversas variantes do VRP, modelando diferentes si\-tu\-a\-ções. Situações distintas trazem características e requisitos também distintos. A variante mais estudada assume uma capacidade nos veículos, podendo ser um peso ou volume máximo que um veículo pode carregar. Tal restrição implica em um número máximo de entregas que um determinado veículo pode realizar. A capacidade pode ser variável ou igual para os veículos. Esta variante é nomeada \Sigla{\textit{Capacitated Vehicle Routing Problem}}{CVRP}.

Assumindo janelas temporais, em que entregas possuem um tempo específico para serem realizadas, temos o \Sigla{\textit{Vehicle Routing Problem with Time Windows}}{VRPTW}.

Também há a possibilidade de utilizar diversos depósitos, estes problemas são chamados de \Sigla{\textit{Multi Depot Vehicle Routing Problem}}{MDVRP}. De acordo com o trabalho \cite{eksioglu2009vehicle} cerca de $11\%$ dos artigos abordados assumem a variante com múltiplos depósitos, contra $90,\!5\%$ do CVRP e $37,\!9\%$ do VRPTW.

O VRP é um problema da classe NP-difícil, isto significa que não é possível obter soluções ótimas em tempo de execução polinomial ao menos que $P=NP$ \cite{talbi2009metaheuristics}.

No problema proposto neste trabalho modelamos um conjunto de entregas de refeições e um conjunto de restaurantes que preparam marmitas. O cliente pede uma marmita de um conjunto limitado de opções, e a refeição é preparada por um restaurante credenciado. O usuário não sabe qual restaurante irá preparar a sua refeição. 

Este é um modelo encontrado na prática em aplicativos como \emph{Ifood}. Tal aplicativo tem uma opção chamada \emph{Loop} que se enquadra no modelo apresentado. Nele, o cliente faz sua requisição no dia anterior à entrega. Uma vez que o \emph{Ifood} fecha o conjunto de entregas, é necessário realizar um processo de otimização em três etapas: divisão das refeições entre os restaurantes; cada restaurante divide suas refeições entre os entregadores; e, para cada entregador, é necessário otimizar a rota de entrega. Todas as três etapas são cobertas por nosso modelo.


\subsection{Modelos de Otimização}
Segundo \textcite{talbi2009metaheuristics}, cientistas, engenheiros e gerentes sempre tomam decisões. Conforme o mundo se torna cada vez mais complexo, o processo de decisão tem de ser feito de forma racional e otimizada.

Para \textcite{talbi2009metaheuristics} o processo de decisão consiste em quatro passos: \textit{(i)} formular o problema; \textit{(ii)} modelar o problema; \textit{(iii)} otimizar o problema; e \textit{(iv)} implementar uma solução.

No primeiro passo, o problema é identificado e é feita uma formulação inicial. Esta formulação pode ser imprecisa. Fatores internos, externos e os objetivos são delineados.

Na segunda etapa, um modelo matemático é construído para o problema, podendo ser inspirado por modelos similares presentes na literatura. Nestes casos, reduz-se o problema para modelos mais bem estudados. Geralmente todo modelo é uma simplificação da realidade, sendo por vezes, incompleto. A modelagem pode apresentar aproximações e alguns fenômenos podem não ser representados por serem muito complexos ou pouco relevantes para os objetivos pretendidos.

Uma vez que o problema é modelado, deve-se trabalhar em uma boa solução para o problema, podendo ser ótima ou sub-ótima. Após este passo, o processo de decisão pode finalmente ser implementado e a solução encontrada testada. Entretanto, \textit{benchmarks} permitem que o desempenho dos processos possam ser testados antes de serem implementados.

\subsection{Objetivo e contribuições}

O foco deste trabalho consiste em modelar uma situação realística de entregas de comida. Inicialmente encontramos um problema que consiste em entregas de marmitas na cidade de Rio Claro-SP. Cada marmita é requisita de um aplicativo. Após este passo, as marmitas são divididas entre diversos restaurantes (em torno de 7). Cada restaurante divide as entregas entre motoboys.

O processo de otimização visa minimizar o comprimento das rotas, minimizar o número de entregadores e maximizar o balanceamento das cargas dos restaurantes.

Neste trabalho não demos enfase a algoritmos de otimização para o VRP. Fizemos apenas as buscas locais 2-opt e 3-opt. O foco do trabalho foi o desenvolvimento do \textit{benchmark}.

Inicialmente, partimos de um mapa da cidade de Rio Claro e modelamos cada rua como uma cadeia poligonal. Manualmente, extraímos um total de 1566 ruas. O motoboy anda entre dois pontos do mapa por meio de caminhos mínimos no grafo das ruas.

Uma vez estabelecido o mapa de Rio Claro num grafo planar 2D, as entregas foram geradas aleatoriamente. Na nossa modelagem, não foi atribuído sentido nas ruas. Assume-se que todas as ruas são de mão dupla.

Os depósitos também foram gerados aleatoriamente. Após isso, foi feita a modelagem da solução como um vetor de permutação com entregas e depósitos. Cada permutação representa uma solução completa, com rotas associadas a depósitos. Após definida a representação da solução, foi feito o cálculo das três funções objetivo para cada solução.

O objetivo deste TCC foi criar um \textit{benchmark} baseado nas ruas de Rio Claro, para o problema Multi-Depot VRP. Após as instâncias serem lançadas na internet, outros alunos irão desenvolver técnicas de otimização para o problema.

Este trabalho foi uma extensão do trabalho \textcite{2020zeni}, em que o VRP para entregas dos carteiros na cidade de Artur Nogueira foi modelado, pois abrange a metodologia desenvolvida nele e a amplia para gerar o VRP contendo múltiplos depósitos.

\section{Levantamento bibliográfico}\label{chp:levantamento}

Ainda antes de surgir o nome VRP, o trabalho \cite{wren1972computer} já apresentou um método para resolver o MDVRP. Segundo esse trabalho, considerando todos os depósitos juntos, geralmente se obtém soluções melhores dividindo o problema em sub-problemas e resolvendo-os separadamente. Também de acordo com esse trabalho, o MDVRP veio ao conhecimento em 1966 em conexão com o problema de posicionamento de depósito. Nesse caso os depósitos considerados eram pequenos, com com cada depósito contendo apenas um veículo.

O artigo \cite{pooley1994integrated} modela um caso real acontecido nos anos 80 sobre uma empresa canadense do ramo de alimentos chamada "Ault Foods". Neste trabalho, o autor expõe desafios na elaboração da estratégia de uma mudança de instalações do setor de fluidos dessa empresa. Segundo o autor, os gestores queriam, entre outros fatores, saber se mudando as instalações permitiriam alcançar uma redução de custos significativa. Para isso, a Ault reuniu uma equipe para estudar esse problema. A equipe tinha dois objetivos: \textit{(i)} analisar as oportunidades de uma racionalização da produção existente e da rede de distribuição e \textit{(ii)} determinar quais novas oportunidades estariam associadas com as novas instalações de produção e distribuição. 

O trabalho é um bom exemplo de MDVRP, pois entre as questões que deveriam ser respondidas estava qual depósito deve servir cada cliente. Após a pesquisa de mudança nas instalações de produção, a equipe identificou que poderia trabalhar com sete depósitos sem perder qualidade no atendimento dos seus clientes, cinco depósitos a menos que no modelo anterior da empresa.

O artigo \cite{cassidy1972tramp} trabalha com a formulação de um problema de entrega de merendas escolares para escolas no Interior de Londres. Os autores descrevem o problema considerando veículos diferentes, com restrições de tempo e de tamanho dos veículos. Esse trabalho também aborda o conceito de \textit{"carried through"}, ou seja, as entregas podem ser recebidas por diferentes depósitos antes de serem entregues.

O trabalho \cite{renaud1996tabu} destaca que o MDVRP consiste na construção de um conjunto de rotas da seguinte forma: \textit{(i)} cada rota começa e termina no mesmo depósito; \textit{(ii)} cada cliente deve ser visitado exatamente uma vez por um veículo; \textit{(iii)} a demanda total de cada rota não pode exceder a capacidade do veículo; (iv) a duração total de cada rota não pode exceder o limite definido e \textit{(v)} a rota total é minimizada. O trabalho também destaca que o MDVRP é um problema muito difícil de encontrar soluções ótimas mesmo para instâncias relativamente pequenas e propõe um algoritmo de busca tabu para este problema.


Em \cite{bowerman1995multi} os autores desenvolvem um novo problema chamado \Sigla{\textit{School Bus Routing Problem}}{SBRP} baseado no VRP porém com a adição de diversos outros objetivos e restrições. O trabalho apresenta cinco objetivos que devem ser minimizados: \textit{(i)} número de rotas; \textit{(ii)} comprimento da rota total; \textit{(iii)} balanceamento do carregamento dos alunos transportados em cada rota; \textit{(iv)} balanceamento do comprimento das rotas e \textit{(v)} distância de caminhada dos estudantes. O problema também apresenta três principais restrições: \textit{(i)} número máximo de alunos por cada rota; \textit{(ii)} uma distância máxima por cada rota e \textit{(iii)} um tempo máximo que os estudantes podem caminhar.

\textit{Benchmarks} são ferramentas muito presentes na computação e em outras áreas e permitem avaliar o desempenho de produtos e métodos. \textit{Benchmarks} para o VRP permitem a avaliação de desempenho de métodos e há na literatura diversos \textit{benchmarks} conhecidos. O principal e mais conhecido \textit{benchmark} para o VRP é a TSPLIB, documentada no trabalho \cite{reinelt1991tsplib}. Esse trabalho traz um agrupamento de diversas instâncias de diferentes trabalhos, sendo disponibilizadas em \cite{reinelt1995tsplib95}. A TSPLIB possui instâncias para o TSP simétrico e assimétrico, bem como instâncias de problemas relacionados, tais como o CVRP, o problema da ordenação sequêncial e o ciclo de Hamilton.

Segundo \cite{uchoa2017new} na literatura de métodos exatos, tornou-se comum seguir a convenção TSPLIB de arredondamento de distâncias para o número inteiro mais próximo (Reinelt, 1991) e também de fixar o número de rotas.


Em \cite{meira2020assess} há o desenvolvimento de um \textit{benchmark} para o VRP. O trabalho aborda o problema de entrega de correspondências por carteiros na cidade de Artur Nogueira, denominado PostVRP. O PostVRP considera veículos sem capacidade e uma rota máxima limitada a um tempo de 6 a 8 horas e apenas um depósito. Esse trabalho desenvolve uma metodologia de geração probabilística de clientes em um espaço 2d e um método para o cálculo de distâncias. O \textit{benchmark} gerado contém instâncias com até 30 mil clientes. No PostVRP há uma abordagem multi-objetivo, pois considera três objetivos: \textit{(i)} a minimização do comprimento da rota total; \textit{(ii)} a minimização do número de veículos e (\textit{iii}) a redução da variabilidade no comprimento das rotas.

Em \cite{karakativc2015survey} os autores fazem uma pesquisa sobre algoritmos genéticos para resolver o problema do MDVRP. Os autores em \cite{karakativc2015survey} destacam que os trabalhos \cite{arostegui2006empirical} e \cite{youssef2001evolutionary} comparam meta-heurísticas de busca tabu, simulated annealing e algoritmos genéticos para o \textit{facility location problem} e \textit{location planning}. Nos trabalhos destacados os algoritmos genéticos apresentam vantagens de tempo e desempenho em comparação com outras meta-heurísticas em situações com restrições de tempo e limitação de poder computacional, ainda que muitas vezes outras meta-heurísticas são capazes de trazer soluções melhores.

\section{Metodologia}\label{chp:metodologia}
Neste capítulo, descreveremos o desenvolvimento da pesquisa. Este trabalho é uma extensão de \cite{2020zeni} pois utiliza a metodologia utilizada nesse trabalho para ampliá-la a um novo problema, com adição de novas características. Na Seção \ref{defVRP} definimos a variante FoodDeliveryVRP. Na Seção \ref{modelCreate} abordaremos a criação do modelo. A Seção \ref{benchCreate} aborda aspectos da criação do \textit{benchmark}. Na Seção \ref{solIni} abordaremos soluções iniciais. Não demos o foco deste trabalho para algoritmos e otimizações e sim para o desenvolvimento do \textit{benchmark}.

\subsection{Definindo a variante}\label{defVRP}
Considere um grafo ponderado completo $G(V,E)$ e uma função de custo $w : V \times V \rightarrow \mathbb{N}$. O conjunto dos depósitos é dado por $\Pi \in V$. 
Em nosso problema, cada depósito representa um restaurante. Suponha que existem 10 restaurantes e 100 entregas de marmita. Podemos associar cada entrega a um restaurante de maneira livre. Os restaurantes produzem marmitas idênticas.

O conjunto de clientes é dado por $C=V\setminus\Pi$ e seu número, por $n=|C|$. O conjunto de clientes é representado por $C=\lbrace c_1,\ldots,c_n\rbrace$ e o conjunto de depósitos é representado por $\Pi =\{\pi_1,\pi_2,\ldots,\pi_{n'}\} $. Podemos definir o conjunto de vértices como $$V=\{v_1,\ldots,v_{n+n'}\}=\{\pi_1,\ldots,\pi_{n'},c_1,\dots,c_n\},$$
ou seja, os primeiros $n'$ vértices são depósitos e os $n$ vértices seguintes são pontos de entrega.

Existe um valor $k_{max} \in \mathbb{N}$ que representa o número máximo de veículos para cada depósito. Existe também um valor $R_{max}$ que é o comprimento máximo da rota. Rotas mais longas que $R_{max}$ são consideradas inviáveis.

Uma instância do problema é representada por um grafo completo $G(V,E,w)$, por dois inteiros $n$ e $n'$, por um inteiro $k_{max}$ e por um inteiro $R_{max}$
Os depósitos são os primeiros $n'$ elementos de $V$ e os pontos de entrega são os $n$ vértices seguintes.

A função $w : V \times V \rightarrow \mathbb{N}$ representa o custo entre qualquer par de elementos de $V$. Está implícito neste grafo que o custo $w$ é obtido por meio do algoritmo de caminhos mínimos sobre as ruas de Rio Claro. Neste momento, o processamento já foi feito e temos apenas o custo $w$. O valor de $w$ é dado em unidades de tempo. Já está embutido o tempo de deslocamento mais o tempo de entrega.

Considere uma sequência $$S(C,k_{max}) = (c_1, \ldots, c_n,\pi_1, \ldots, \pi_1,\pi_2, \ldots, \pi_2,\ldots,
\pi_{n'}, \ldots, \pi_{n'})$$  montada através da inserção de todos os elementos de $C$ em $S(C,k)$.  Após isso cada vértice em $\Pi$ é inserido  $k_{max}$ vezes. Cada permutação de $S(C,k)$ representa uma solução do VRP.

A entrega que está na posição $s_i\in S(C,k_{max})$ está associada ao primeiro depósito a sua esquerda, caso exista, ou $\pi_0$ caso contrário. 
Cada rota $R_i$ é obtida por uma sequência de entregas entre o $i$-ésimo depósito e o $(i+1)$-ésimo depósito (ou entre o $i$-ésimo depósito e o fim da sequencia). Toda rota $R_i$ começa no depósito $\pi_i$ associado a ela.

Por exemplo, considere a solução $$S' = (c_1,c_2,c_3,c_4,\pi_a,c_5,c_6,c_7,\pi_b,c_8,c_9,c_{10}).$$ Nesse exemplo, $R_1 = (\pi_0,c_1,c_2,c_3,c_4,\pi_0)$, $R_2 = (\pi_a,c_5,c_6,c_7,\pi_a)$ e $R_3 = (\pi_b,c_8,c_9,c_{10},\pi_b)$. Seja $Particao(S)=(R_1,\ldots,R_{k'})$ o conjunto gerado pela quebra da sequencia original em rotas. A quebra se dá toda vez que um depósito é encontrado na sequencia original. Por definição rotas sem clientes não fazem parte de $Particao(S)$. Ou seja, $Particao(c_1,c_2,\pi_a,\pi_b,c_3,c_4)$ é $$\lbrace(\pi_0,c_1,c_2,\pi_0),(\pi_b,c_3,c_4,\pi_b)\rbrace$$ e não $\lbrace(\pi_0,c_1,c_2,\pi_0),(\pi_a,\pi_a),(\pi_b,c_3,c_4,\pi_b)\rbrace$.
O tamanho de uma rota $R = (r_1, \ldots, r_m)$, em que $m$ representa a posição do último vértice presente na rota, é dado por:

\begin{center}
    $W(R) =  \displaystyle\sum_{i=1}^{m-1}w(r_i, r_{i+1})$.
\end{center}

O comprimento de uma solução $S = (s_1, \ldots, s_m)$ é calculado como a soma do comprimento das rotas, como a seguir:

\begin{center}
    $W(S) = \displaystyle\sum_{R \in Particao(S)} W(R)$
\end{center}

Depois de feita a partição das rotas, cada entrega está associada a um depósito.
Seja $Entregas(\pi_i)$ as entregas associadas ao depósito $\pi_i\in \Pi$. O número de clientes de $\pi_i$ é dado por:

\begin{center}
     $|Entregas(\pi_i)|$
\end{center}

O número de veículos usados em uma dada solução é igual ao número de rotas não vazias, ou seja $|Particao(S)|$. 

Definimos três funções de custo $f_1(S)$, $f_2(S)$ e $f_3(S)$. Em $f_1(S)$ utilizamos o comprimento da solução, com $f_1(S) = W(S)$. 

Na função $f_2(S)$ estabelecemos o número de veículos com $$f_2(S) = |Particao(S)|.$$ 

Por fim, a função $f_3(S)$ mede o grau de variabilidade entre a quantidade de entregas de cada restaurante, a fim de distribuir de maneira mais uniforme a carga de trabalho dos restaurantes. Para esta finalidade, seu resultado se dá através do calculo do desvio padrão:
\begin{center}
$f_3(S) = \displaystyle\sqrt{\dfrac{\displaystyle\sum_{\forall \pi_i} (|Entregas(\pi_i)| - \overline{|Entregas|})^2}{|\Pi|}}$
\end{center}

A seguir, temos a definição formal do problema:

\begin{definition}{\textbf{FoodDeliveryVRP}.}
Dado um conjunto de elementos $V$, uma função de custo $w :  V \times V \rightarrow \mathbb{N} $,  uma constante $ k_{max} \in S $ representa o número máximo de veículos por depósito, o comprimento máximo da rota $R_{max} \in \mathbb{N}$ e dois inteiros $n'$ e $n$. Seja $V=\{v_1,\ldots,v_{n+n'}\}=\{\pi_1,\ldots,\pi_{n'},c_1,\dots,c_n\}$. Considere a sequência $S(C,k_{max})$ e $Pe$  o conjunto de todas as permutações factíveis de $S(C,k_{max})$ respeitando o $R_{max}$. O FoodDeliveryVRP consiste na minimização de $(f_1(S'), f_2(S'), f_3(S')) $ para todo $S'\in Pe$.
\end{definition}

\subsection{Criação do Modelo}\label{modelCreate}

Assim como no trabalho \cite{2020zeni}, as ruas são modeladas como uma cadeia poligonal $P$ definida como um conjunto de coordenadas planares, tal que $P = (c_1, \ldots, c_n)$, sendo $c \in \mathbb{R}^2$ para todo $c \in P$. O grafo $G(V,E)$ é gerado a partir do conjunto dessas ruas. Cada vértice $v \in V$ está associado a uma coordenada cartesiana $(x,y) \in \mathbb{R}^2$ e cada aresta $e = (u,v)$ é um segmento de reta entre $u$ e $v$. O custo da aresta é definido por distância euclidiana.

Diferentemente do trabalho \cite{2020zeni}, não atribuímos custo para atravessar a rua, tendo em vista que as entregas são realizadas por veículo e esse custo seria desprezível. Entretanto, utilizamos a atribuição de densidade desenvolvida no trabalho \cite{2020zeni} para atribuir probabilidades de ruas receberem entregas. Este método é usado para que ruas centrais tenham maior  probabilidade de receberem entregas comparadas com ruas isoladas por unidade de comprimento.

\subsection{Criação do Benchmark}\label{benchCreate}

O trabalho de \textcite{2020zeni} apresenta uma ferramenta para criação de \textit{benchmarks}, essa ferramenta possui três arquivos de configuração:
\begin{itemize}
    \item \textit{Background.png}: contém uma imagem utilizada para a melhor visualização. Essa imagem também serve de base para a construção do modelo;
    \item \textit{Model.txt}: contém informações sobre o modelo, tais como: custo adicional para realizar uma entrega, valores de conversão de pixel, atributos das ruas, local do depósito e o mapa das ruas;
    \item \textit{Instances.txt}: contém as configurações de cada instância, tais como número máximo de veículos, quantidade de entregas e tamanho máximo de cada rota.
\end{itemize}

A criação do \textit{benchmark} ocorreu com base em modificações da ferramenta desenvolvida por \cite{2020zeni}. As modificações ocorreram nos arquivos de configuração: \textit{background.png}; \textit{model.txt} e \textit{instances.txt} além de modificações na própria ferramenta. Outra modificação foi no programa \cite{postvrpbench} que faz o \textit{parsing} do arquivo de instância. Com base nesse programa modificado foram desenvolvidas as soluções propostas na Seção \ref{solIni}. As mudanças foram feitas para que a ferramenta lidasse com múltiplos depósitos.

\begin{figure}[htb!]
    \centering
    \includegraphics[width=\textwidth]{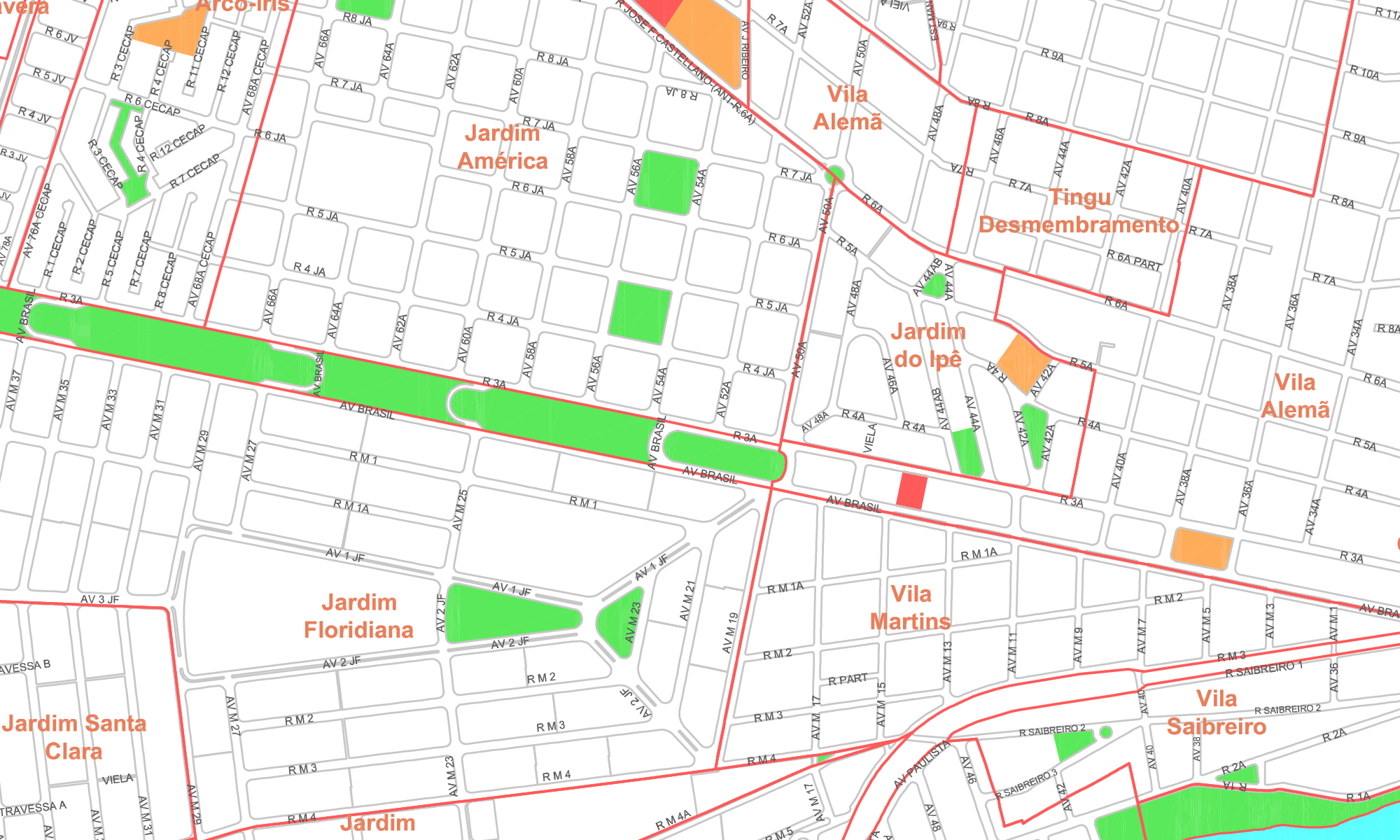}
    \hrule
    \includegraphics[width=\textwidth]{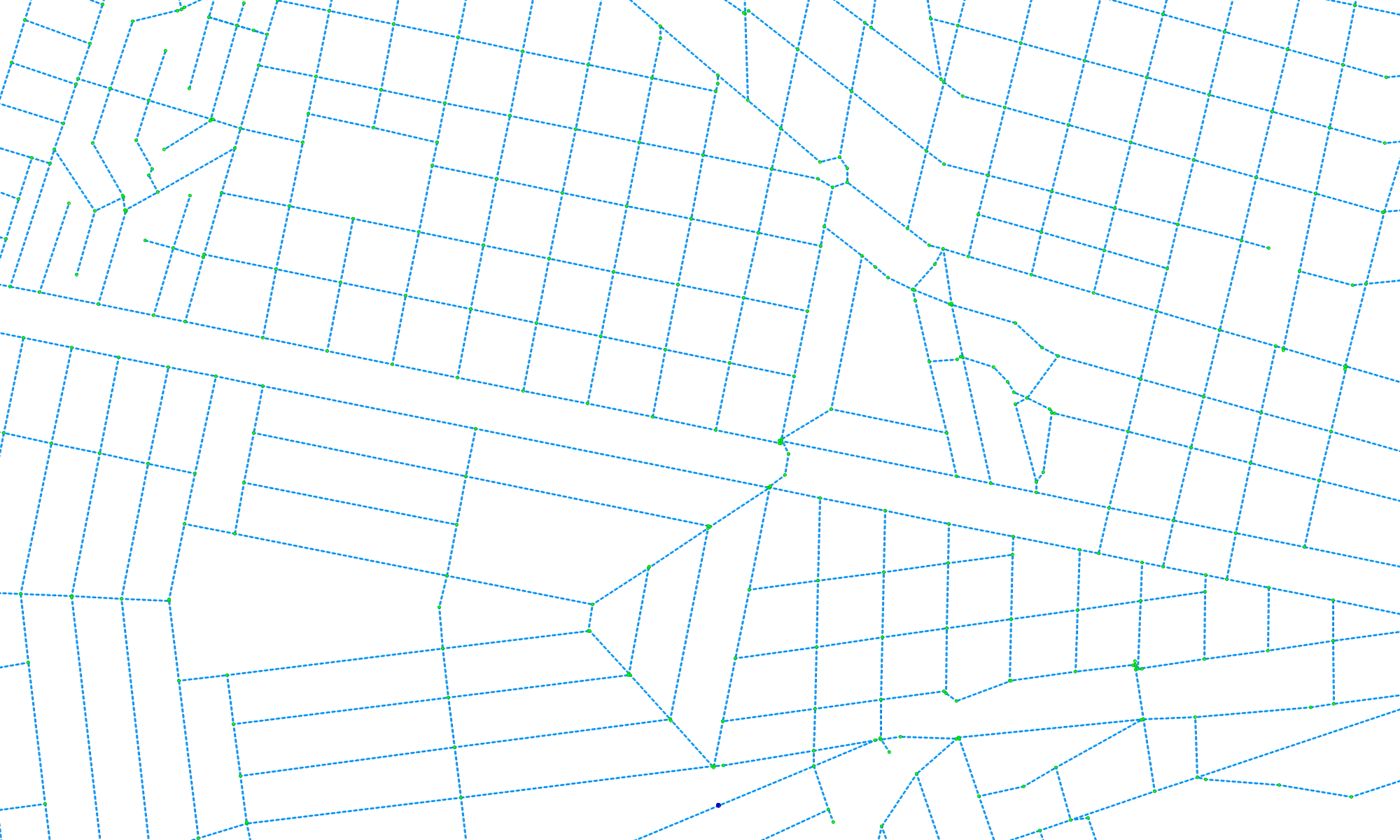}
    \caption{Modelagem das ruas. Acima um recorte do mapa original e abaixo as ruas extraídas (tracejados azuis). Fonte: Prefeitura de Rio Claro (cima) e Próprio Autor (baixo).}
    \label{modelDev}
\end{figure}

Iniciamos o desenvolvimento através da extração manual das coordenadas correspondentes às ruas do mapa da cidade de Rio Claro. Nesta etapa, houve a retirada de 1566 ruas. As coordenadas foram retiradas com o auxílio de um programa de desenho, entretanto para facilitar a extração das coordenadas, desenvolvemos uma ferramenta simples em \textit{javascript}. A ferramenta copia automaticamente as coordenadas do ponto clicado para a área de transferência. Foram desconsiderados na etapa da extração condomínios e atalhos não acessíveis por veículos. A Figura \ref{modelDev} ilustra a extração das ruas e a representação de algumas ruas modeladas através desse processo.

Os depósitos foram obtidos aleatoriamente através do algoritmo de geração de entregas, portanto no arquivo \textit{model.txt} não apresenta valores de coordenadas do depósito, ao invés disso, o arquivo \textit{instances.txt} apresenta o atributo \textit{number of depots} que determina o número de depósitos referente a cada instância. As primeiras $|\Pi|$ entregas geradas são interpretadas como depósitos (ou seja, restaurantes).

A construção do modelo foi com base em um recorte de uma imagem digital do mapa de Rio Claro-SP disponibilizada no site da prefeitura\footnote{Mapa da cidade de Rio Claro disponível em: \url{https://www.rioclaro.sp.gov.br/municipio/publicacoes_mapas.php}}.


Da mesma forma do trabalho \cite{2020zeni} cada rua foi classificada utilizando os atributos Região (R), Tipo (T) e Zona (Z). A classificação das ruas ocorreu de forma arbitrária\footnote{Como arbitrária entenda o conhecimento do aluno sobre a cidade de Rio Claro. O aluno mora na cidade de Santa Gertrudes, vizinha a Rio Claro.}, tendo em vista a quantidade de ruas ser muito grande para uma verificação manual. As probabilidades atribuídas, penalidades e nomenclaturas utilizadas também seguem o padrão definido por \cite{2020zeni} conforme a mostra a Tabela \ref{tab:ruasAtributos}.

\begin{table}[!htp]
\caption[Atributos das ruas]{Atributo, nível e valores das penalidades. Fonte: \cite{2020zeni}.}
\label{tab:ruasAtributos}
\begin{center}
\begin{tabular}{ccccc}
\hline 
Atributo & Nível 1$_{(pen)}$ & Nível 2$_{(pen)}$ & Nível 3$_{(pen)}$ & Nível 4$_{(pen)}$\\ \hline 
Região(R) & central$_{(1.0)}$ & periférico$_{(0.75)}$ & distante$_{(0.4)}$ & isolado$_{(0.2)}$ \\
Tipo(T) &  avenida$_{(1.0)}$ & rua$_{(0.75)}$ & alameda$_{(0.4)}$ & rodovia$_{(0)}$\\
Zona(Z) & comercial$_{(1.0)}$ & misto$_{(0.75)}$ & residencial$_{(0.4)}$ & - \\\hline 
\end{tabular}
\end{center}
\end{table}

Cada rua de Rio Claro recebeu um atributo Região, dentre os valores \emph{central} (penalização 1), \emph{periférico} (penalização 0.75), \emph{distante} (penalização 0.4) e \emph{isolado} (penalização 0.2). A penalização 0.4 significa que a densidade de entregas é multiplicada pelo valor 0.4. 

Cada rua de Rio Clario recebeu um atributo \emph{Tipo}, com valores \emph{avenida}, \emph{rua}, \emph{alameda} e \emph{rodovia}, com penalizações $(1;0.75;0.4;0)$ respectivamente.

Finalmente, cada rua de Rio Claro recebeu um atributo \emph{Zona} com valores \emph{comercial}, \emph{misto} e \emph{residencial} com penalizações $(1;0.75;0.4)$ respectivamente.

\subsection{Soluções iniciais}\label{solIni}
Neste trabalho não demos foco para a otimização e busca por soluções, entretanto implementamos uma solução inicial com a simples atribuição de entregas ao depósito mais próximo. Também houve a implementação dos algoritmos 2-opt e 3-opt. Estes são algoritmos de busca local conhecidos na literatura e que serão abordados mais detalhadamente na Seção \ref{sec:localBusca}.

\subsubsection{Busca local}\label{sec:localBusca}


Segundo \cite{laporte2009fifty}, as meta-heurísticas de busca local exploram o espaço de soluções através do movimento da solução atual para outra solução na vizinhança. Para \cite{laporte2009fifty}, os passos mais importantes de uma estratégia de busca local é determinar as regras que definem o conceito de vizinhança e o mecanismo para explorá-la. Para Talbi \cite{talbi2009metaheuristics} a busca local começa de uma dada solução inicial e, a cada iteração, a heurística substitui a solução atual por uma vizinha que melhora a função objetivo. A busca local para quando todos os candidatos vizinhos são piores que a solução atual, ou seja, o ótimo local foi atingido. 

\begin{figure}[htb!]
    \centering
    \includegraphics[scale=0.8]{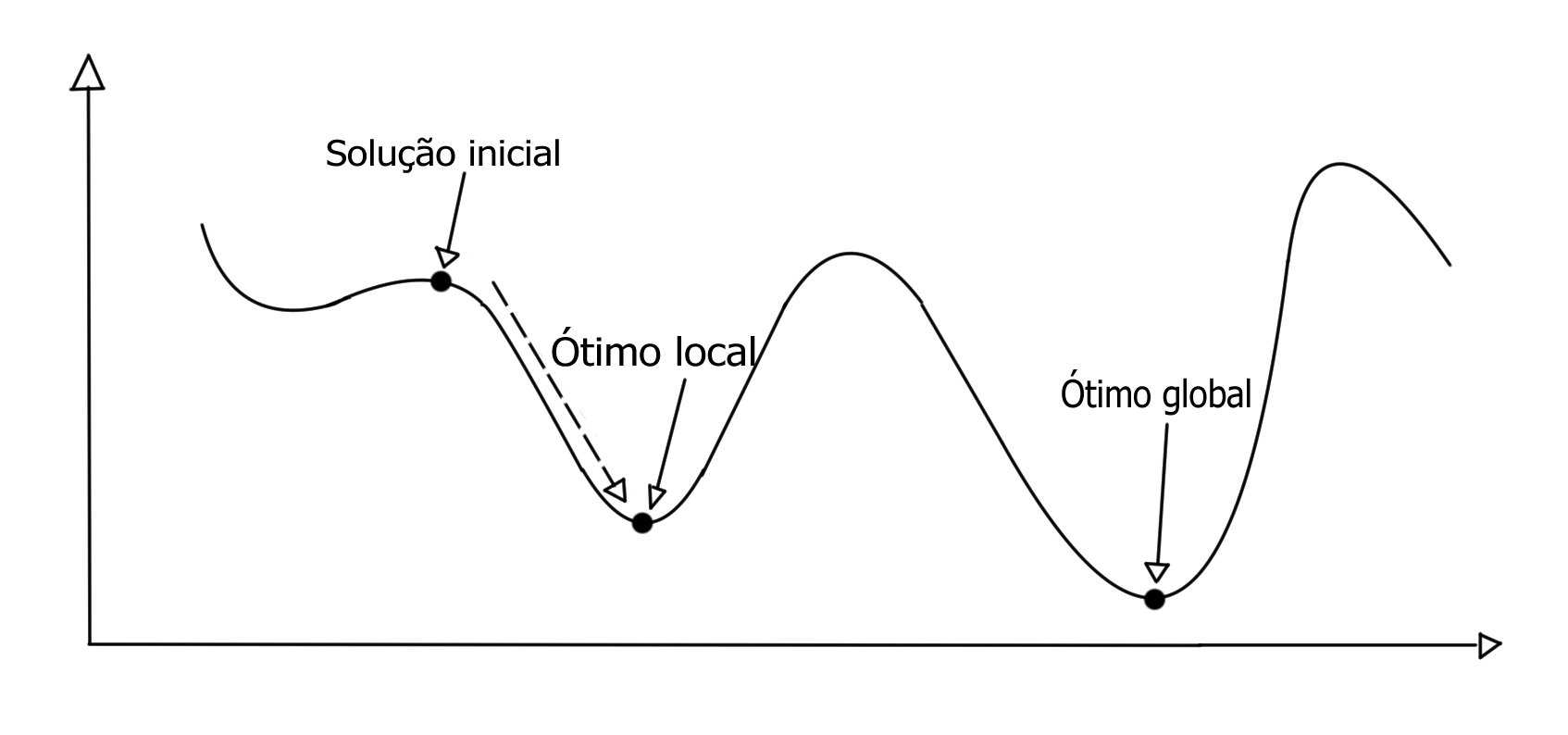}
    \caption{Busca local: ótimo local $\times$ ótimo global. Fonte: Próprio Autor.}
    \label{buscalocal}
\end{figure}

Para determinar uma vizinhança em uma busca local, vários fatores devem ser considerados, pois, ao passo em que uma vizinhança mais abrangente pode gerar soluções melhores, também trazem um maior custo computacional. Outro importante fator a ser considerado é como as melhorias serão escolhidas, Talbi \cite{talbi2009metaheuristics} descreve três possíveis estratégias para a seleção:
\begin{itemize}
    \item \textbf{Melhor melhora:} Nesta estratégia, é escolhida a melhor solução vizinha possível, através da avaliação todos os possíveis movimentos de troca. 
    \item \textbf{Primeira melhora:} Nesta estratégia, é escolhida a primeira solução encontrada melhor que a solução atual e então realiza a troca.
    \item \textbf{Escolha aleatória:} Consiste na escolha aleatória aplicada entre as soluções vizinhas. É feita a troca quando apresentam uma melhora com relação a solução atual.
\end{itemize}

Talbi \cite{talbi2009metaheuristics} cita que as estratégias de primeira melhora apresentam uma melhor relação entre qualidade das soluções e tempo de busca quando aplicadas  a soluções geradas aleatoriamente enquanto que as estratégias de melhor melhora são melhores para soluções geradas através de um método guloso.

Neste trabalho utilizaremos a estratégia da melhor melhora.

\subsubsection{2-OPT}
\begin{figure}
\centering
\includegraphics[scale=1.0]{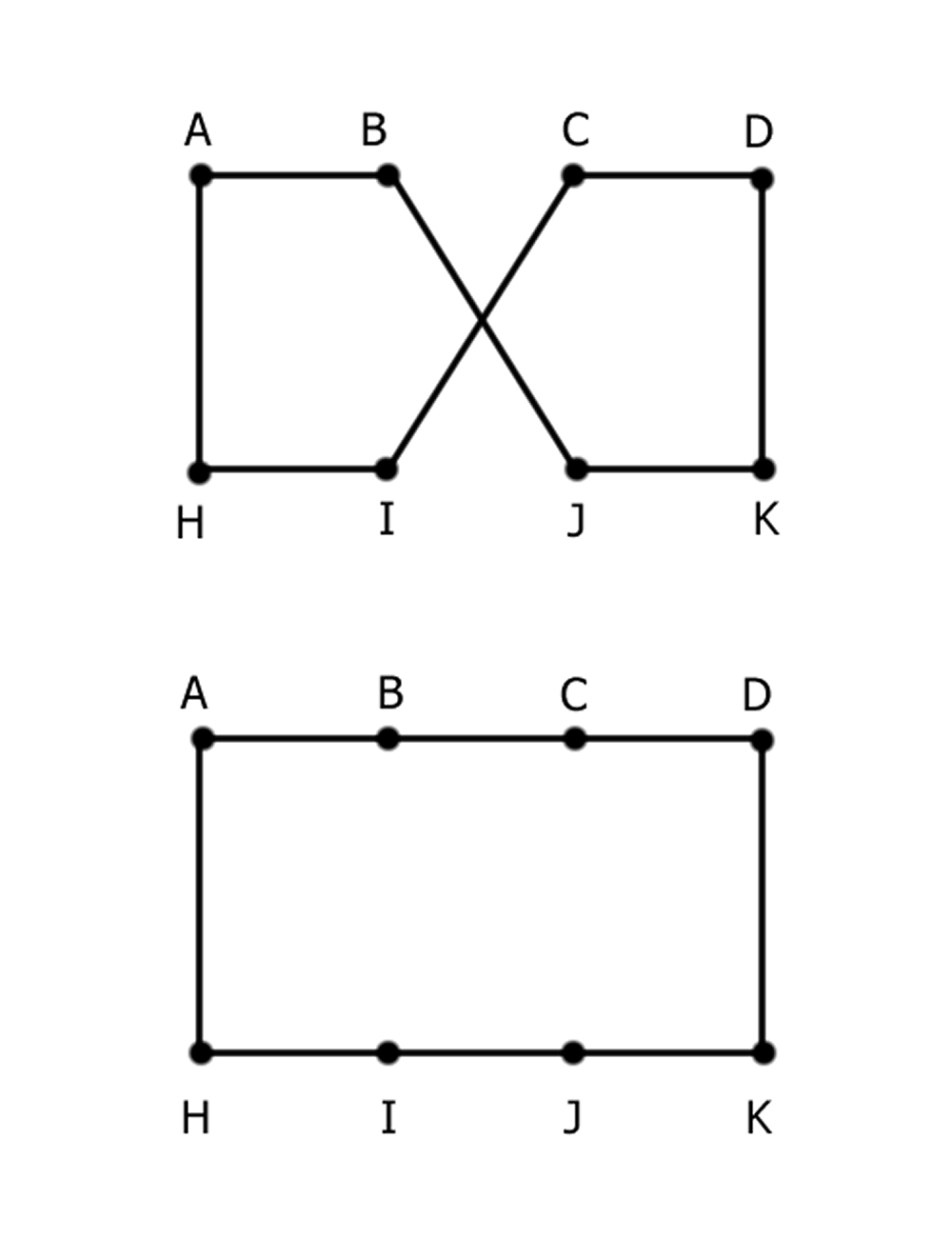}
\caption{Vizinhança do algoritmo 2-OPT. Fonte: Próprio Autor.}
\label{opt2swap}
\end{figure}

O algoritmo 2-opt (algoritmo \ref{2optAlg}) foi inicialmente proposto em 1958 por Croes \cite{croes1958method}, mas a forma básica do movimento foi sugerida por Flood em 1956 \cite{flood1956traveling}, ambos para o TSP. O 2-opt é um algoritmo de busca local que tem como objetivo a remoção de cruzamentos entre arestas. Seu conceito de vizinhança é representado pela troca de duas arestas, ou seja, dada uma solução $R$, $R'$ é considerada vizinha de $R$ se $R'$ for possível de ser obtida a partir da troca de duas arestas em $R$ como na Figura \ref{opt2swap}.

\begin{algorithm}[htb!]
\SetAlgoLined
\DontPrintSemicolon 
\KwIn{Uma rota $R$ com $n$ elementos}
\KwOut{Uma rota $R'$ com o comprimento minimizado}
$melhora \gets verdadeiro$\\
inicio\\
\While{melhora é verdadeiro}{
    \For{$R \in S$}{
        \For{$a \gets 0$ \textbf{até} $n-2$}{
            \For{$b \gets a+2$ \textbf{até} $n$}{
                $delta \gets$ diferença de distância entre as antigas arestas e as novas\\
                \If{delta \text{$<$} 0}{
                    realiza a troca\\
                    $melhora \gets verdadeiro$\\
                }
            }
        }
    }
    
    $melhora \gets falso$\\
}
\caption{2-opt}
\label{2optAlg}
\end{algorithm}

\subsubsection{3-OPT}
\begin{figure}[htb!]
\centering
\includegraphics[width=0.8\textwidth]{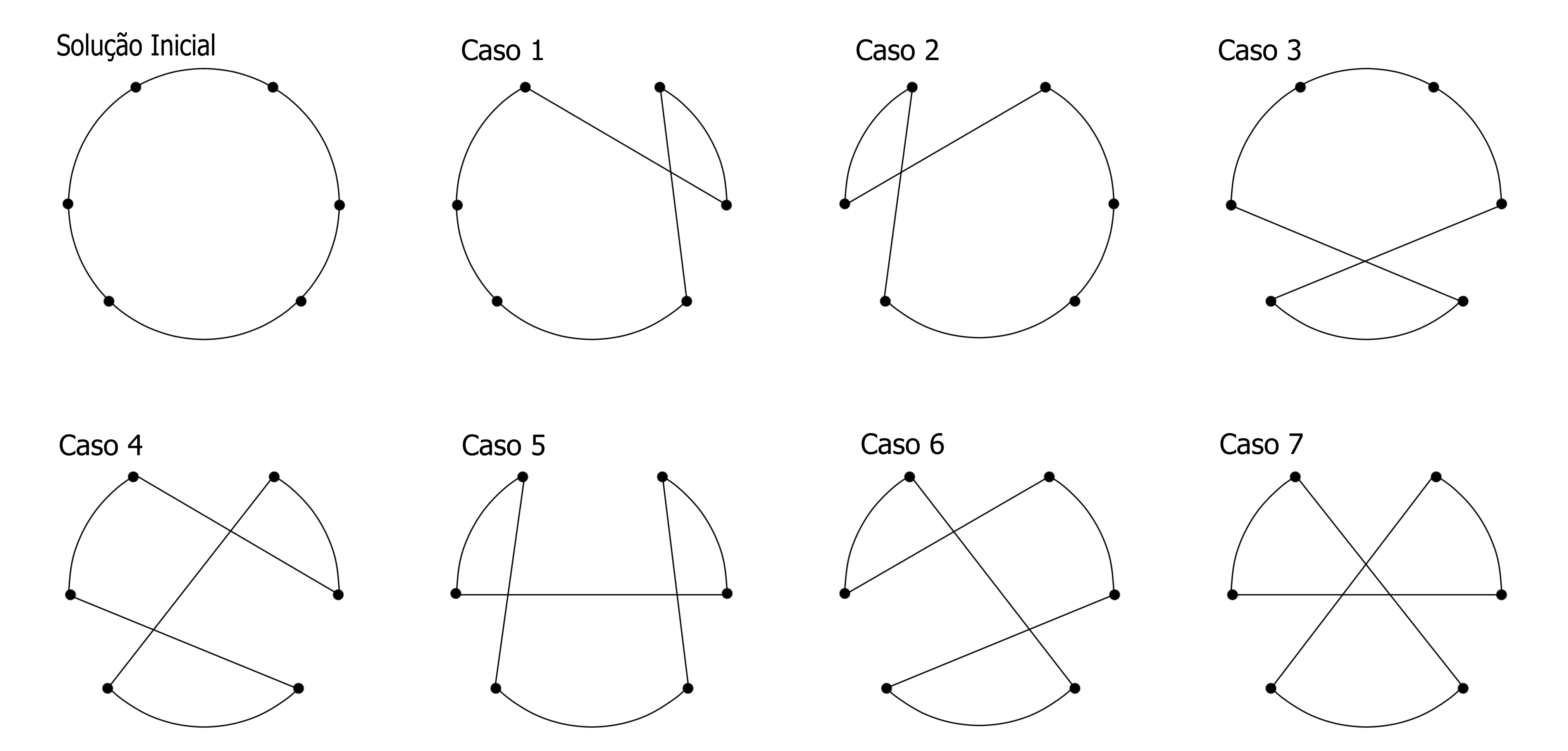}
\caption{Vizinhança do algoritmo 3-OPT. Fonte: Próprio Autor.}
\label{opt3swap}
\end{figure}

O algoritmo 3-opt é uma generalização do algoritmo 2-opt através da expansão do conceito de vizinhança de duas arestas para três. O algoritmo foi descrito pela primeira vez em uma palestra apresentada na 14ª ORSA de 1958 por F. Bock \cite{johnson1997traveling}, sendo mais tarde descrito e utilizado no trabalho de Lin em 1965 \cite{lin1965computer}. Para o 3-opt, uma solução $R'$ é vizinha de $R$ se ela pode ser obtida a partir de $R$ pela troca de até três arestas. Por causa disso, neste algoritmo, todos os possíveis casos de troca entre três arestas são compreendidos. A Figura \ref{opt3swap} representa todos os casos de trocas gerando soluções vizinhas de uma dada solução inicial. Considerando três arestas, temos um total de sete trocas possíveis.

\begin{algorithm}[htb!]
\SetAlgoLined
\DontPrintSemicolon 
\KwIn{Uma rota $R$ com $n$ elementos}
\KwOut{Uma rota $R'$ com o comprimento minimizado}
\KwData{$melhora \gets verdadeiro$}
inicio\\
\While{melhora é verdadeiro}{
    \For{$a \gets 0$ \textbf{até} $n-4$}{
      \For{$b \gets a+2$ \textbf{até} $n-2$}{
        \For{$c \gets b+2$ \textbf{até} $n$}{
            encontrar o melhor caso de troca\\
            $delta \gets$ diferença de distância entre as antigas arestas e as novas\\
            \If{delta \text{$<$} 0}{
                realiza a troca\\
                $melhora \gets verdadeiro$\\
            }
        }
      }
    }
    $melhora \gets falso$\\
}
\caption{3-opt}
\label{3optAlg}
\end{algorithm}

\section{Resultados}\label{chp:resultados}

Neste capitulo será apresentado o \textit{benchmark} desenvolvido neste trabalho. Na Seção \ref{benchResults} o modelo será apresentado. Na Seção \ref{benchResults} apresentamos todas as instâncias geradas. Na Seção \ref{sec:solIniResults} discutiremos sobre resultados dos algoritmos 2-opt e 3-opt e a distribuição de entregas.

\subsection{Food Delivery VRP}\label{benchResults}

\begin{figure}[htb!]
    \centering
    \includegraphics[width=\textwidth]{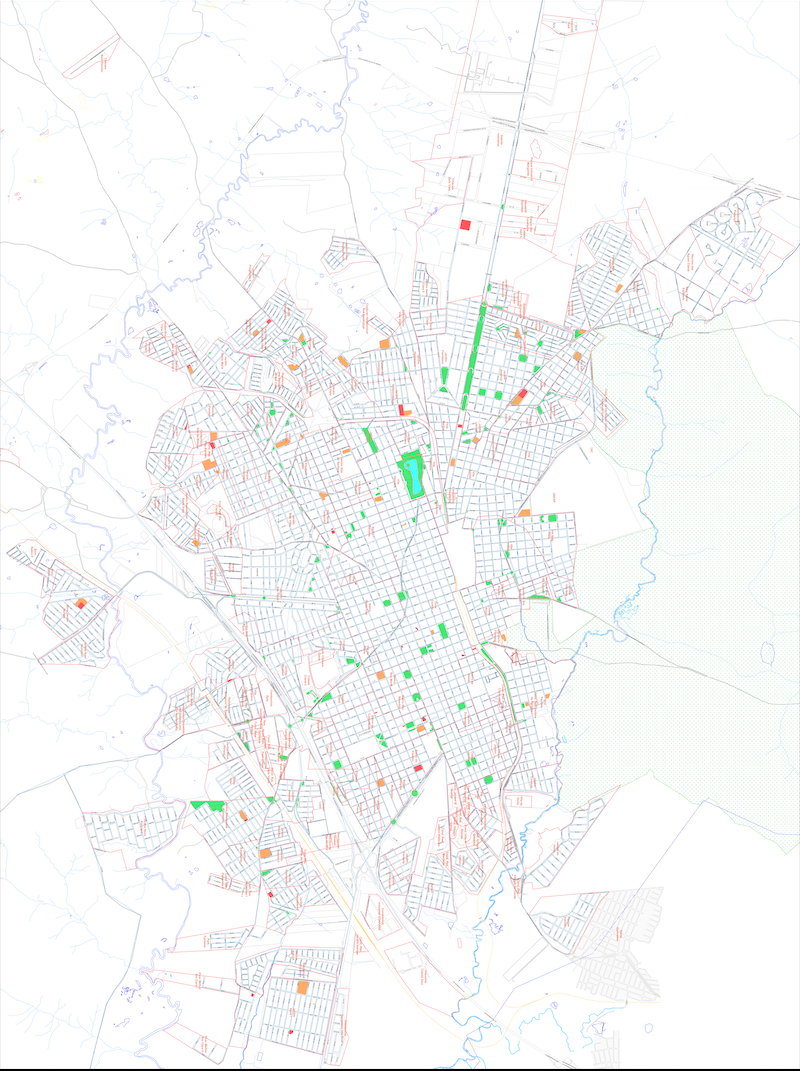}
    \caption{Mapa completo da cidade de Rio Claro. Fonte: Site da prefeitura de Rio Claro.}
    \label{inst23img}
\end{figure}

\begin{figure}[htb!]
    \centering
    \includegraphics[width=\textwidth]{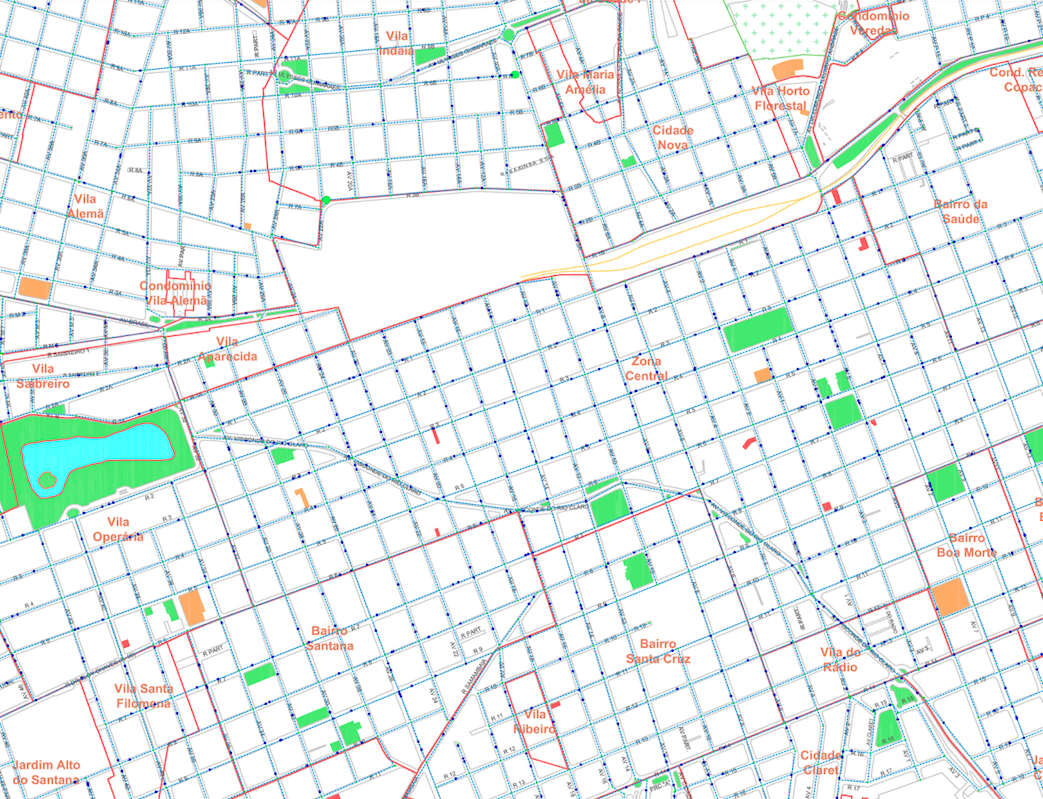}
    \caption{FoodDeliveryVRP instância 23 com 2000 pontos de entrega. Recorte da região central. Cada ponto azul corresponde a uma entrega. Fonte: Próprio Autor.}
    \label{instancia23centro}
\end{figure}

\begin{figure}[htb!]
    \centering
    \includegraphics[width=0.9\textwidth]{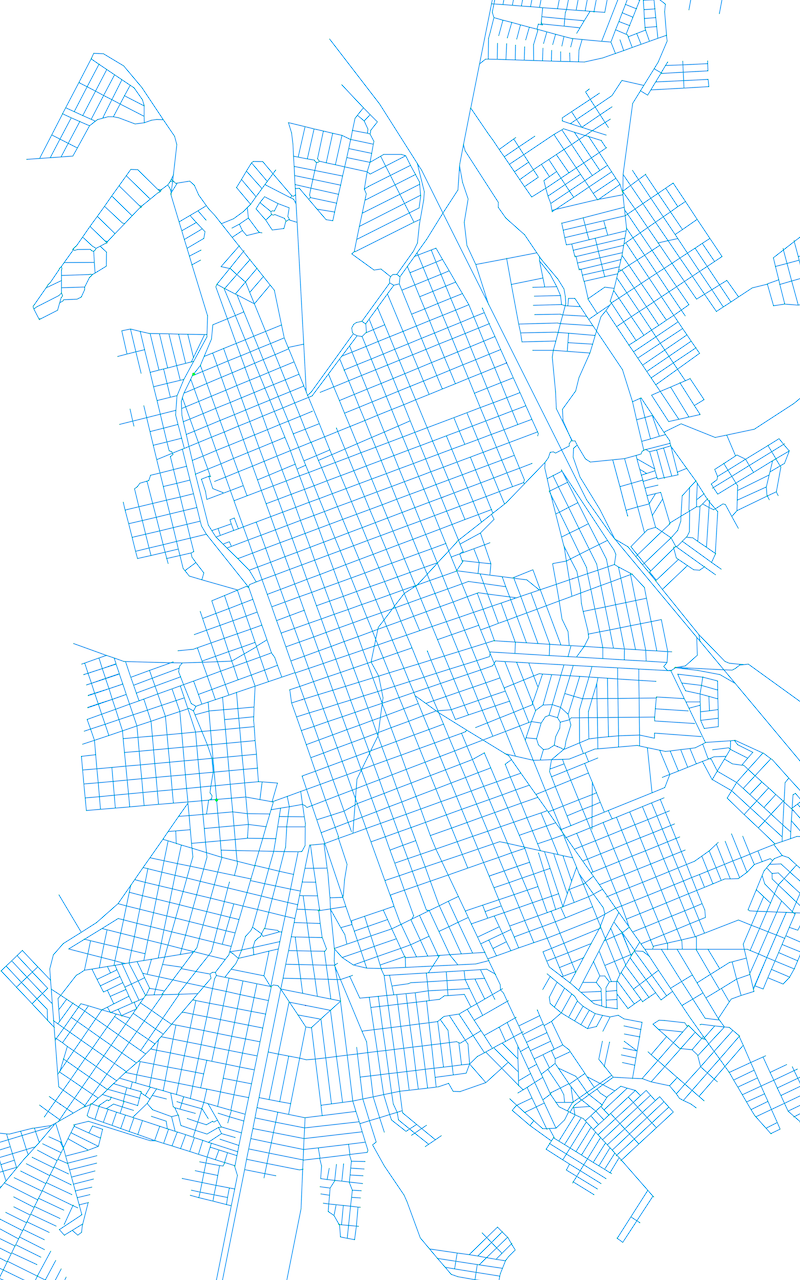}
    \caption{Ruas do modelo. Fonte: Próprio Autor.}
    \label{blankmap}
\end{figure}

O FoodDeliveryVRP contém 23 instâncias iniciando com 10 pontos de entrega e 2 depósitos e terminando em 2000 pontos de entrega e 7 depósitos. O número de depósitos é variável para diferentes instâncias, mesmo as que contém o mesmo número de entregas. A Figura \ref{inst23img} apresenta uma imagem completa do modelo gerado do \textit{benchmark}. Na Figura \ref{instancia23centro} podemos observar um recorte da região central para uma instância com 2000 pontos de entrega. A Figura \ref{blankmap} apresenta o mapa completo apenas com as ruas modeladas representadas pelas linhas tracejadas em azul.

\begin{figure}[htb!]
    \centering
    \includegraphics[width=\textwidth]{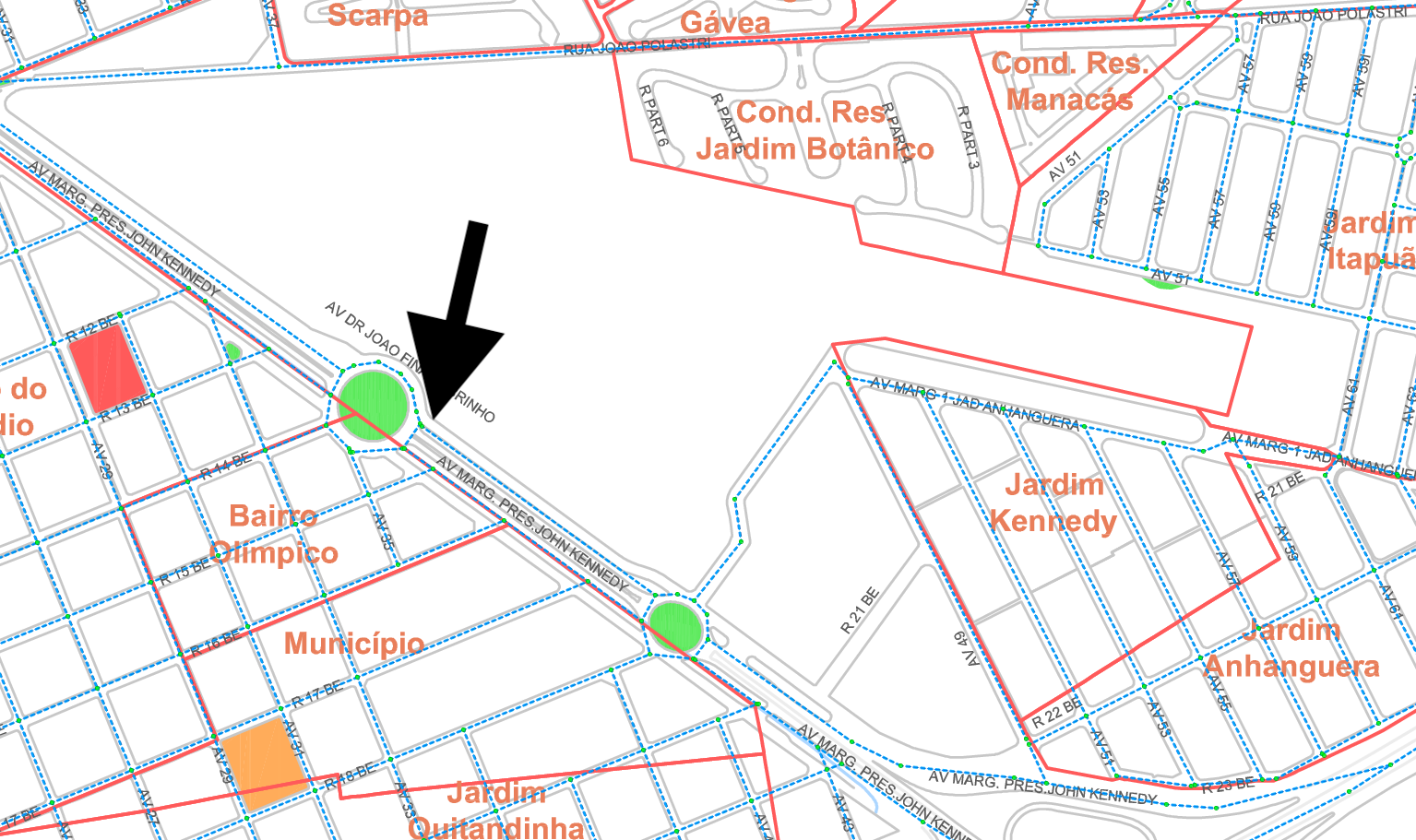}
    \caption{Avenida Marginal Presidente John Kennedy. A seta aponta a rua descrita pelas coordenadas. Fonte: Próprio Autor.}
    \label{ruaExemplo}
\end{figure}

Cada rua é modelada como uma sequência de coordenadas para representar segmentos de reta. Ruas com várias curvas ou circulares são modeladas através da inserção de vários vértices ao longo de seu percurso. Como exemplo, a Avenida Marginal Presidente John Kennedy (parte de cima) no mapa da Figura \ref{ruaExemplo} é modelada através da sequência:

AV. PRES JOHN KENNEDY [CENTRAL,LARGEAVENUE,COMMERCIAL]
[4062,8629]-[4086,8652]-[4366,9040]-[4360,9052]-[4356,9079]-[4369,9103]-[4386,9115]-[4410,9120]-[4426,9126]-[4609,9379]-[4608,9403]-[4615,9421]-[4630,9435]-[4657,9436]-[4768,9592]-[4794,9645]-[4807,9768]-[4836,9933]-[4840,9987]-[4794,10173]

\begin{figure}[htb!]
    \centering
    \includegraphics[width=\textwidth]{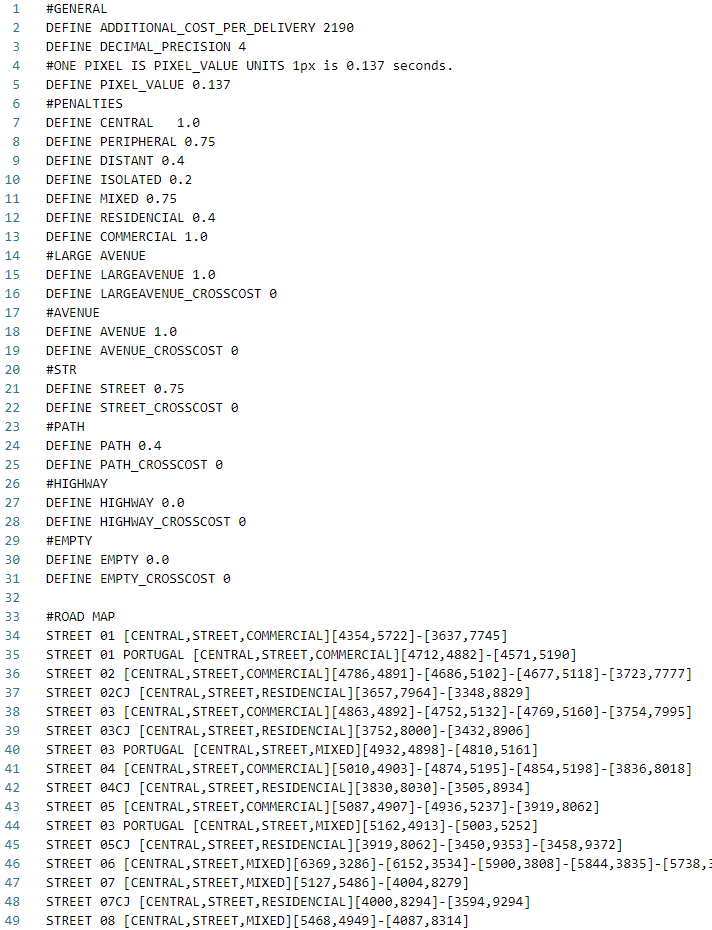}
    \caption{Inicio do arquivo \textit{model.txt}. Fonte: Próprio Autor.}
    \label{modelINI}
\end{figure}

\begin{figure}[htb!]
    \centering
    \includegraphics[width=\textwidth]{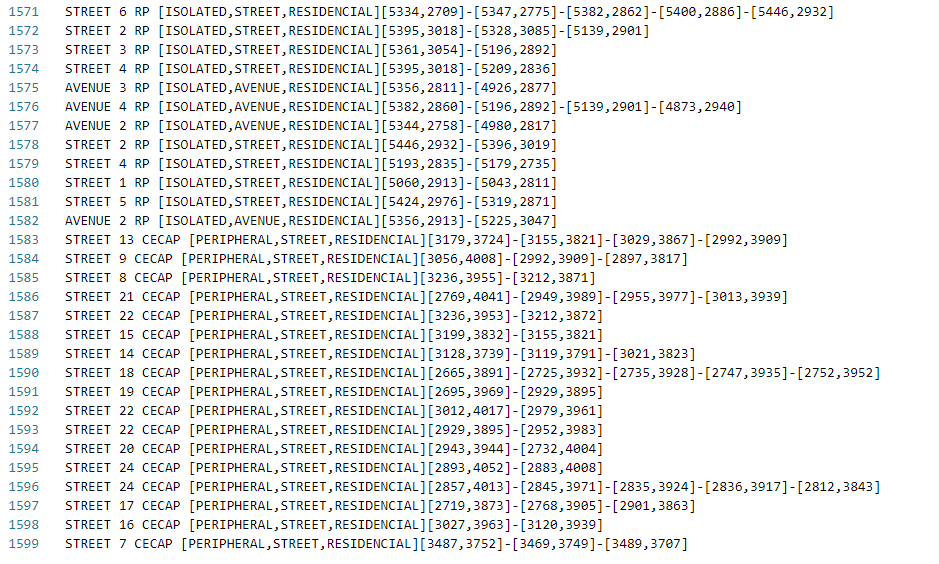}
    \caption{Fim do arquivo \textit{model.txt}. Fonte: Próprio Autor.}
    \label{modelFIM}
\end{figure}

O arquivo final \emph{model.txt} com os parâmetros para gerar as instâncias estão descritos nas Figuras~\ref{modelINI} e \ref{modelFIM}. Na Figura~\ref{modelINI} temos as primeiras 49 linhas do arquivo e na Figura~\ref{modelFIM} temos as últimas linhas do arquivo, terminando na linha 1599.

A velocidade dos veículos ($PIXEL\_VALUE$) foi definida no arquivo \emph{model.txt} em $0.137s/pixel$ representando uma motocicleta a $30Km/h$. 
Esta constante serve para relacionar diretamente pixel com segundo.
O custo adicional por entrega foi definido em 2190 pixeis, o equivalente a 5 minutos por entrega. A rota máxima permitida foi de 13138 pixeis, o equivalente a meia hora. As refeições no mundo real possuem um tempo muito curto para serem entregues. Fatores como a perda de temperatura podem afetar a qualidade da refeição, impedindo a realização de rotas longas.

Outro fator importante a ser considerado é que a entrega pelo \textit{motoboy} não é tão rápida quanto a entrega de uma correspondência. Para o carteiro realizar uma entrega, na maioria das vezes, deixa a correspondência na caixa de correio, enquanto que o \emph{motoboy} precisa ser atendido pelo cliente e também cobrar pelo serviço, aumentando muito o custo de cada entrega. Todos esses fatores tornam o FoodDeliveryVRP um problema com grandes restrições de rota, tornando inviáveis instâncias muito grandes, como será mostrado na Seção \ref{sec:solIniResults}.

O mapa gerado pode ser alterado para gerar novas instâncias outras variantes do VRP como coleta de lixo ou entrega dos carteiros. Por exemplo, podemos usar o mapa para gerar novas instâncias para o PostVRP \cite{2020zeni}. Para isso, no arquivo \textit{model.txt} deve haver a adição de um depósito fixo, ajuste no valor de pixel para corresponder a velocidade de um carteiro e alteração da precisão decimal. Também é necessário gerar um novo arquivo \textit{instances.txt}.


A Tabela \ref{tab:instancias1} apresenta as instâncias do FoodDeliveryVRP. São um total de 23 instâncias com até 2000 pontos de entrega e até 7 depósitos.

\begin{table}[!htp]
\caption[Instâncias do Benchmark]{Instâncias do FoodDeliveyVRP.}
\label{tab:instancias1}
\begin{center}
\begin{tabular}{ccccc}
\hline 
ID & Instância & Entregas & Depósitos & Veículos(max)\\ \hline 
1 & FoodDelivery\_10\_0 & 10 & 2 & 5 \\
2 & FoodDelivery\_10\_1 & 10 & 3 & 5 \\
3 & FoodDelivery\_10\_2 & 10 & 4 & 5\\
4 & FoodDelivery\_20\_0 & 20 & 2 & 5\\
5 & FoodDelivery\_20\_1 & 20 & 3 & 5\\
6 & FoodDelivery\_20\_2 & 20 & 4 & 5\\
7 & FoodDelivery\_50\_0 & 50 & 4 & 5\\
8 & FoodDelivery\_50\_1 & 50 & 5 & 5\\
9 & FoodDelivery\_50\_2 & 50 & 6 & 5\\
10 & FoodDelivery\_100\_0 & 100 & 5 & 10\\
11 & FoodDelivery\_100\_1 & 100 & 6 & 10\\
12 & FoodDelivery\_100\_2 & 100 & 7 & 10\\
13 & FoodDelivery\_200\_0 & 200 & 5 & 10\\
14 & FoodDelivery\_200\_1 & 200 & 6 & 10\\
15 & FoodDelivery\_200\_2 & 200 & 7 & 10\\
16 & FoodDelivery\_500\_0 & 500 & 5 & 20\\
17 & FoodDelivery\_500\_1 & 500 & 6 & 20\\
18 & FoodDelivery\_500\_2 & 500 & 7 & 20\\
19 & FoodDelivery\_1000\_0 & 1000 & 5 & 40\\
20 & FoodDelivery\_1000\_1 & 1000 & 6 & 40\\
21 & FoodDelivery\_1000\_2 & 1000 & 7 & 40\\
22 & FoodDelivery\_2000\_0 & 2000 & 7 & 80\\
23 & FoodDelivery\_2000\_1 & 2000 & 7 & 80\\

\hline 
\end{tabular}
\end{center}
\end{table}

\subsection{Resultados Iniciais}\label{sec:solIniResults}

As soluções iniciais foram obtidas através da utilização de um algoritmo guloso que atribui as entregas ao depósito mais próximo. Esta etapa será chamada de \emph{clusterização}.
Note que no problema descrito, uma entrega pode ser feita por qualquer restaurante em $\Pi$.
Os algoritmos de otimização atuam sobre a rota resultante da clusterização para gerar soluções vizinhas iterativamente até o ótimo local ser alcançado. O processo não distribui entregas em veículos, há apenas a divisão em depósitos. As tabelas \ref{tab:instancia23results}, \ref{tab:instancia17results}, \ref{tab:instancia13results} e \ref{tab:instancia7results} apresentam os resultados obtidos das instâncias 23, 17, 13 e 7 respectivamente. Os resultados de tempo foram obtidos através da média de 5 execuções em um processador de 3,4GHz.

Devido ao comprimento máximo ser de 30 minutos e o tempo de entrega de 5 minutos, o FoodDeliveryVRP apresenta problemas com instâncias com muitas entregas, pois na prática cada rota pode atender um número máximo 5 clientes, tornando o número de rotas necessárias bastante elevado.

\begin{table}[!htp]
\caption[Resultados da instância 23]{Resultados da instância 23 com 2000 clientes e 7 depósitos. Os algoritmos 2-opt e 3-opt otimizam a rota após a clusterização. Cada restaurante possui uma única rota que não foi quebrada em veículos. Os valores de $f_2$ e $f_3$ obtidos pela clusterização foram de 7 e 213,29 respectivamente.}
\label{tab:instancia23results}
\begin{center}
\begin{tabular}{cccc}
\hline 
Algoritmo & $f_1$ (h) & Tempo Médio & Desvio Padrão\\
 & Soma Rotas  & Execução (ms) & Tempo Execução (ms) \\\hline
Cluster. & 298,9 & $<1$ & $<1$\\
Cluster. $+$ 2-opt & 178,6 & 336 & 17\\
Cluster. $+$ 3-opt & 177,8 & 47590 & 8185\\
\hline 
\end{tabular}
\end{center}
\end{table}

\begin{table}[!htp]
\caption[Resultados da instância 17]{
Resultados da instância 17 com 500 clientes e 6 depósitos. Os algoritmos 2-opt e 3-opt otimizam a rota após a clusterização. Cada restaurante possui uma única rota que não foi quebrada em veículos. Os valores de $f_2$ e $f_3$ obtidos pela clusterização foram de 6 e 71,45 respectivamente.
}
\label{tab:instancia17results}
\begin{center}
\begin{tabular}{cccc}
\hline 
Algoritmo & $f_1$ (h) & Tempo Médio & Desvio Padrão\\
 & Soma Rotas  & Execução (ms) & Tempo Execução (ms) \\\hline
Cluster. & 77,1 & $<1$ & $<1$\\
Cluster. $+$ 2-opt & 48,5 & 47 & 6,7\\
Cluster. $+$ 3-opt & 48,2 & 891 & 103\\
\hline 
\end{tabular}
\end{center}
\end{table}

\begin{table}[!htp]
\caption[Resultados da instância 13]{
Resultados da instância 13 com 200 clientes e 5 depósitos. Os algoritmos 2-opt e 3-opt otimizam a rota após a clusterização. Cada restaurante possui uma única rota que não foi quebrada em veículos. Os valores de $f_2$ e $f_3$ obtidos pela clusterização foram de 5 e 23,82 respectivamente.}
\label{tab:instancia13results}
\begin{center}
\begin{tabular}{cccc}
\hline 
Algoritmo & $f_1$ (h) & Tempo Médio & Desvio Padrão\\
 & Soma Rotas  & Execução (ms) & Tempo Execução (ms) \\\hline
Cluster. & 30,6 & $<1$ & $<1$\\
Cluster. $+$ 2-opt & 21,5 & 15 & $<1$\\
Cluster. $+$ 3-opt & 21,2 & 156 & 8,41\\
\hline 
\end{tabular}
\end{center}
\end{table}

\begin{table}[!htp]
\caption[Resultados da instância 7]{
Resultados da instância 7 com 50 clientes e 4 depósitos. Os algoritmos 2-opt e 3-opt otimizam a rota após a clusterização. Cada restaurante possui uma única rota que não foi quebrada em veículos. Os valores de $f_2$ e $f_3$ obtidos pela clusterização foram de 4 e 7,78 respectivamente.
}
\label{tab:instancia7results}
\begin{center}
\begin{tabular}{cccc}
\hline 
Algoritmo & $f_1$ (h) & Tempo Médio & Desvio Padrão\\
 & Soma Rotas  & Execução (ms) & Tempo Execução (ms) \\\hline
Cluster. & 8,6 & $<1$ & $<1$\\
Cluster. $+$ 2-opt & 6,8 & $<1$ & $<1$\\
Cluster. $+$ 3-opt & 6,7 & 15 & 8,41\\
\hline 
\end{tabular}
\end{center}
\end{table}

\begin{figure}[htb!]
    \centering
    \includegraphics[width=0.9\textwidth]{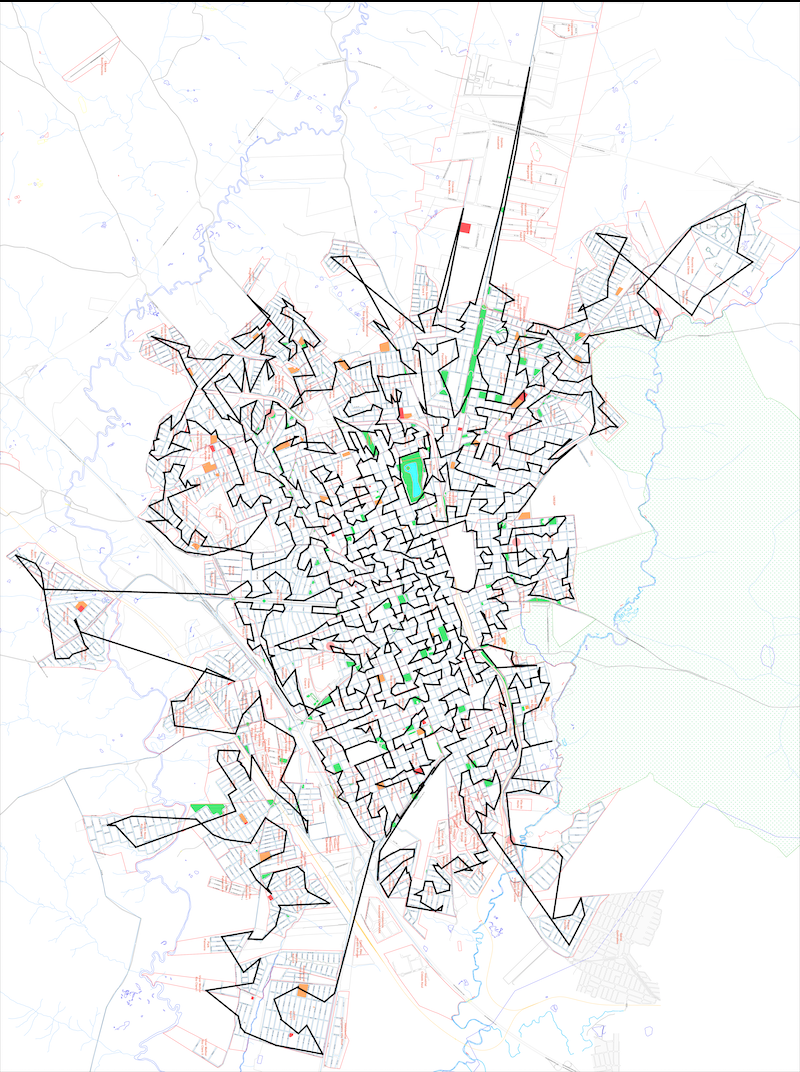}
    \caption{Resultado dos algoritmos de divisão de rota seguido do 3-opt para a instância 23. Foi feita uma clusterização em sete grupos, um para cada depósito e depois foi executado o $3$-opt em cada grupo. Fonte: Próprio Autor.}
    \label{inst23w3opt}
\end{figure}

A Figura \ref{inst23w3opt} representa o resultado obtido a partir da otimização utilizando a clusterização e o algoritmo 3-opt.

\section{Conclusões}\label{chp:conclusoes}





Problemas de roteamento são problemas de logística. Neste trabalho modelamos um \textit{benchmark} para a variante MDVRP, estendendo o trabalho PostVRP \cite{2020zeni} no qual ocorre a modelagem de um problema envolvendo a entrega de correspondências por carteiros através de um único depósito na cidade de Artur Nogueira-SP.

O FoodDeliveryVRP trabalha com três objetivos: \emph{(i)} o comprimento total da rota, \emph{(ii)} o número de veículos (rotas) e \emph{(iii)} a variação da quantidade de entregas entre depósitos. Ao final foi gerado um modelo contendo 1566 ruas. Foram desconsiderados sentidos das ruas, assumindo assim que todas as ruas são de mão dupla.

Para este trabalho consideramos um comprimento máximo da rota de 30 minutos. Além disso, o tempo adicional por entrega é de 5 minutos. Sendo assim, o FoodDeliveryVRP é um problema com grande número de veículos, pois  cada veículo pode fazer no máximo 5 entregas. 

O FoodDeliveryVRP  pode ser estendido no futuro através da adição de janelas temporais ou entregas dinâmicas, podendo também servir de \textit{benchmark} para o desenvolvimento de algoritmos de otimização para outros problemas.

O foco deste trabalho não foi o desenvolvimento de algoritmos, entretanto foram implementados os algoritmos de busca local conhecidos na literatura 2-opt e 3-opt e um algoritmo guloso de clusterização.

Analisando os resultados de soluções iniciais através da divisão por depósitos e otimização através dos algoritmos 2-opt e 3-opt, podemos observar que a diferença nos resultados entre os algoritmos 2-opt e 3-opt foram de menos de 1,1\% sendo o 2-opt o algoritmo mais rápido enquanto que o 3-opt apresentou melhores resultados. 

\section{Agradecimentos}

Esta pesquisa foi parcialmente financiada pela FAPESP (Auxílio  2015/11937-9) e CAPES (Código Financiador 001).

\bibliographystyle{plain}

\bibliography{bibliografia}

\end{document}